\sffamily
\documentclass[11pt,twoside]{article}

\usepackage{algorithmic}       
\usepackage{amsbsy}            
\usepackage{amsmath}           
\usepackage{amssymb}           
\usepackage{bm}                
\usepackage{color}             
\usepackage[T1]{fontenc}       
\usepackage{fancyvrb,moreverb} 
\usepackage{float}             
\usepackage{gloss}             
  \setglosslabel{\bfseries#2\ifglossshort{ (#3)}{}} 
  \renewcommand\glossheading[1]{\stopglosslist}     
     \makegloss   
\usepackage{graphicx}          
\usepackage{hyphenat}          
\usepackage{ifpdf}             
\usepackage{mathrsfs}          
\usepackage{multirow}          
\usepackage{pst-all}           
\usepackage{tikz}              
   \usetikzlibrary{arrows}
   \usetikzlibrary{plothandlers}
   \usetikzlibrary{plotmarks}
   \usetikzlibrary{shapes}
   \usetikzlibrary{backgrounds}
\usepackage{theorem}           
\usepackage{url}               
{\theorembodyfont{\rmfamily} \newtheorem{definition}{Definition}[section]}
{\theorembodyfont{\rmfamily} \newtheorem{example}   {Example}   [section]}
{\theorembodyfont{\rmfamily} \newtheorem{lemma}     {Lemma}     [section]}
{\theorembodyfont{\rmfamily} \newtheorem{remark}    {Remark}    [section]}
{\theorembodyfont{\rmfamily} }
{\theorembodyfont{\rmfamily} \newtheorem{algorithm} {Algorithm} [section]}

\newlength{\textheightDefault    }
\newlength{\textwidthDefault     }
\newlength{\evensidemarginDefault}
\newlength{\oddsidemarginDefault }
\newlength{\topmarginDefault     }
\newlength{\parindentDefault     }
\newlength{\parskipDefault       }

\newcommand{\noopsort}[1]{}  
\def\btrNum{0}
\newcommand{\btr}
  {\ifnum\btrNum=0 Budget Truck Rental (BTR)\def\btrNum{1}\else BTR\fi}
\def\pilNum{0}
\newcommand{\pil}
  {\ifnum\pilNum=0 \textit{Planned Inventory Level} (PIL)\def\pilNum{1}%
     \else PIL\fi}

\providecommand{\def}{definition} %
\DeclareSymbolFont{lettersA}{U}{txmia}{m}{it}
 \DeclareMathSymbol{\natural} {\mathord}{lettersA}{"8E}
 \DeclareMathSymbol{\integer} {\mathord}{lettersA}{"9A}
 \DeclareMathSymbol{\rational}{\mathord}{lettersA}{"91}
 \DeclareMathSymbol{\real}    {\mathord}{lettersA}{"92}
 \DeclareMathSymbol{\complex} {\mathord}{lettersA}{"83}
 \DeclareMathSymbol{\field}   {\mathord}{lettersA}{"86}
 \DeclareMathSymbol{\preal}   {\mathord}{lettersA}{"90}

\newcommand{\milpN}{0}
\newcommand{\milp}{\ifnum\milpN=0%
\mbox{\textit{Mixed-Integer Linear Programming}} (MILP)\else MILP\fi%
\renewcommand{\milpN}{1}}

\newcommand{\slantfrac}[2]{\ensuremath{%
  \kern0.1em\raisebox{0.5ex}{$\scriptstyle #1$}\kern-0.1em
  /\kern-0.15em\raisebox{-0.25ex}{$\scriptstyle #2$}}}

\newcommand{\trademark}{\textsuperscript{\textregistered}}
\newcommand{\excelN}{0}
\newcommand{\excel}{\ifnum\excelN=0\textsc{Excel}\trademark\else\textsc{Excel}\fi%
\renewcommand{\excelN}{1}}

\newcommand{\gamsN}{0}
\newcommand{\gams}{\ifnum\gamsN=0\textsc{gams}\trademark\else\textsc{gams}\fi%
\renewcommand{\gamsN}{1}}
\newcommand{\matlabN}{0}
\newcommand{\matlab}{\ifnum\matlabN=0\textsc{Matlab}\trademark\else\textsc{Matlab}\fi%
\renewcommand{\matlabN}{1}}
\newcommand{\miktex}{\ifnum\miktexN=0\textsc{miktex}\trademark\else\textsc{miktex}\fi%
\renewcommand{\miktexN}{1}}
\newcommand{\pgfN}{0}
\newcommand{\pgf}{\ifnum\pgfN=0\textsc{pgf} (portable graphics format)\else\textsc{pgf}\fi%
\renewcommand{\pgfN}{1}}
\newcommand{\pstN}{0}
\newcommand{\pst}{\ifnum\pstN=0\textsc{PS Tricks}\trademark\else\textsc{PS Tricks}\fi%
\renewcommand{\pstN}{1}}
     
\begin{document}
\title{\LARGE \textbf{Algorithm for solving optimization problems with Interval Valued Probability Measure }}
\author{Phantipa Thipwiwatpotjana\\
        \textit{\normalsize Department of Mathematical Sciences}\\
        \textit{\normalsize University of Colorado Denver}\\
        {\normalsize Phantipa.Thip@cudenver.edu}
        \and
        Weldon A. Lodwick\\
        \textit{\normalsize Department of Mathematical Sciences}\\
        \textit{\normalsize University of Colorado Denver}\\
        {\normalsize Weldon.Lodwick@cudenver.edu} }
\date{}  
\maketitle
\thispagestyle{empty}\begin{abstract}
We are concerned with three types of uncertainties: probabilistic, possibilitistic and interval. By using possibility and necessity measures as an Interval Valued Probability Measure (IVPM), we present IVPM's interval expected values whose possibility distributions are in the form of polynomials. By working with interval expected values of independent uncertainty coefficients in a linear optimization problem together with operations suggested in Lodwick and Jamison \cite{LodJam2007}, the problem after applying these operations becomes a linear programming problem with constant coefficients. This is achieved by the application of two functions. The first is applied to the interval coefficients, $v: I \rightarrow \mathrm{R}^k$, where $I= \left\{\left[a,b\right]\; |\; a \leq b\right\}$. The second is $u:\mathrm{R}^k \rightarrow \mathrm{R}$, applied to the product we got from a previous function.  Similar concepts hold for any types of optimization problems with linear constraints. Moreover, it implied that optimization problems containing all three types of uncertainties in one problem can be solved as ordinary optimization problems.
\end{abstract}

\section{Introduction} \label{intro}
An Interval Valued Probability Measure (IVPM) which is generated from the definitions provided by Weichselberger \cite{Weich2000}, is a tool that gives a partial representation for an unknown probability measure. In this paper we expand the idea in \cite{LodJam2007} of using an IVPM to an optimization problem with uncertain coefficients. The types of these uncertainties in this research are probabilistic, possibilistic and interval uncertainties. We provide the necessary definitions and explanations in the next section. 

To apply in an optimization problem, 
we will use the interval expected value described in section \ref{ivpm} as the representative of each uncertain coefficient. We construct a general form of an interval expected value whose IVPM construction forms a polynomial possibility density function.

In section \ref{lp}, by assuming independence, we apply interval expected values to  all uncertain coefficient random variables. We create an algorithm for a linear program (LP) with interval uncertainty coefficients. Using this algorithm, we show that our uncertainty problem becomes an ordinary LP. Similar details are given in  section \ref{nlp}. Examples, conclusion and further research idea are provided.

\section{Interval Valued Probability Measure (IVPM)} \label{ivpm}
Before giving the definition of an IVPM, we would like the readers to feel comfortable with the notation $\breve{m}$.
\begin{definition}
Define an uncertain random variable $\breve{m}$ as $\breve{m} =\left\{ \begin{array}{l}
\bar{m};\; \mathrm{or}\\
\widehat{m};\; \mathrm{or}\\
\left[m_1,m_2\right]
\end{array}\right.;$
where $\bar{m}$ is corresponding to $m$ as a random variable with 
probabilistic distribution  and $\widehat{m}$ means 
that $m$ is a random variable, with possibilistic distribution. 
\end{definition}

In real applications, we might not know (with certainty) the probability measure for our problems. Lodwick and Jamison \cite{LodJam2007} use an IVPM, $i_{\breve{m}}(A) = \left[i_{\breve{m}}^{-}(A), i_{\breve{m}}^{+}(A)\right]$, to measure a partial representation for an unknown probability measure. 
The original paper for the idea of IVPM is adopted by Weichselberger \cite{Weich2000}. 
We use the following notation and information throughout the paper unless stated otherwise:
\begin{itemize}
\item The arithmetic operations applied to intervals are those of interval arithmetic \cite{Moore1979}. 
\item The set of all intervals contained in $[0,1]$ is denoted as
\[\mathrm{Int}_{[0,1]}\equiv \left\{[a,b] \;|\; 0 \leq a \leq b \leq 1\right\}.\]
\item $\mathcal{S}$ denotes the universal set and $\mathcal{A}$ is a $\sigma$-algebra of $\mathcal{S}$. Note that $\mathcal{S} = \mathrm{R}$.
\end{itemize}
\begin{definition} \label{def:r-probability}
(Weichselberger \cite{Weich2000}) Given measurable space $\left( \mathcal{S},\mathcal{A}%
\right) $, an interval valued function $i_{\breve{m}}:A\subseteq \mathcal{A}%
\rightarrow \mathrm{Int}_{\left[ 0,1\right] }$ is called an \textit{R-probability} if:%
\begin{itemize}
\item $i_{\breve{m}}\left( A\right) =\left[ i_{\breve{m}}^{-}\left( A\right)
,i_{\breve{m}}^{+}\left( A\right) \right] \subseteq \left[ 0,1\right] $,
\item $\exists $ a probability measure, $\Pr $, on $\mathcal{A}$ such
that $\forall A\in \mathcal{A}$, $\Pr \left( A\right) \in i_{\breve{m}}\left(A\right).$
\end{itemize}
\end{definition}

An R-probability from definition \ref{def:r-probability} is an \textit{IVPM} where $i^{-}_{\breve{m}}$ and $i_{\breve{m}}^{+}$ are constructed from a possibility (fuzzy) density function.
\begin{definition} A function $p:\mathcal{S}\rightarrow \left[ 0,1\right] $ is called a \textit{regular possibility distribution function} if 
\[ \sup \left\{ p\left( x\right) \mid x\in \mathcal{S}\right\} =1.\]
Possibility distribution functions (see \cite{WangKlir1992}) define a
possibility measure, $Pos:\mathcal{S}\rightarrow \left[ 0,1\right] $ where 
\[Pos\left( A\right) =\sup \left\{ p\left( x\right) \mid x\in A\right\} \]
and its dual necessity measure is
\[Nec\left( A\right) =1-Pos\left( A^{c}\right) ,\]
where $\;\sup \left\{ p\left( x\right) \mid x\in \emptyset \right\} =0$.  A
necessity distribution function $n:\mathcal{S} \rightarrow \left[ 0,1\right] $ can be
defined by setting 
\[n\left( x\right) =1-p\left( x\right) \]
and the corresponding necessity measure 
\[Nec\left( A\right) =\inf \left\{ n\left( x\right) \mid x\in A^{c}\right\} ,\]
where $\;\inf \left\{ n\left( x\right) \mid x\in \emptyset \right\} =1$. 
\end{definition}
In \cite{LodJam2002}, it is shown that possibility distributions can be constructed
which satisfy the following consistency definition.
\begin{definition} \label{def:consistent}
Let $p:\mathcal{S}\rightarrow \left[ 0,1\right] $ be a regular possibility
distribution function with associated possibility measure $Pos$ and
necessity measure $Nec$. \ Then $p$ is said to be \textit{consistent} with
random variable $X$ if for every measurable set $A$, $Nec\left( A\right)
\leq \Pr \left( X\in A\right) \leq Pos\left( A\right) $.
\end{definition}
The R-probability function $i_{\breve{m}}$ in definition \ref{def:r-probability} is used to define IVPMs. A possibility and necessity pair, $i_{\breve{m}}(A) = \left[Nec(A),Pos(A)\right]$, constructed by definition \ref{def:consistent} is able to  bound an unknown probability of interest. Therefore it can be used to define an IVPM. 
The reader could find more explanations, examples and  a construction of an IVPM in \cite{LodJam2007} .
\subsection{The Interval Expected Value Constructed From Polynomial Possibility Density Function} \label{polyIVPM}
In this paper we consider the interval expected value (definition is given in \cite{LodJam2007}),
\begin{eqnarray}
\int_{\mathrm{R}} x di_{\breve{m}} = \left[\int_{-\infty}^{\infty}xd^{-}(x) dx,\int_{-\infty}^{\infty}xd^{+}(x)dx\right] \label{eq:expected}, 
\end{eqnarray}
 of an IVPM constructed from possibility and necessity measures as the specific upper and lower cumulative probability distribution functions, respectively. When the interval expected value is calculated as we will see, the lower cumulative distribution function gives the right end-point of the interval expected value, while the upper cumulative distribution function defines the left end-point. We give a formal definition of the interval expected value in definition \ref{def:IEV}
\begin{definition} \label{def:IEV}
\textit{The interval expected value} is defined in (\ref{eq:expected}), where  $d^{-}(x) $ refers to  the left possibility density function corresponding to the upper cumulative distribution. Similarly, $d^{+}(x)$ is the right possibility density function corresponding to the lower cumulative distribution.
\end{definition}

This definition will become clear when  we present how to calculate the interval expected value through a polynomial possibility density function.
\begin{definition} \label{def:polyfuzzy}
\textit{A polynomial degree $n$ fuzzy number $a/b/c/d$}, is a random number whose value is fuzzy between $a$ and $d$. The corresponding polynomial fuzzy membership function $f: \mathrm{R}\rightarrow [0,1]$ is defined as
\[f(x) = \left\{\begin{array}{ll}
f_L(x)&;\; \forall x \in [a,b)\\
1 & ;\; \forall x \in [b,c]\\
f_R(x) &;\; \forall x \in (c,d]\\
0 &;\; \mathrm{otherwise},
\end{array}\right.\]
where $f_L$ and $f_R$ are defined in table \ref{tbl:polyDistribution}. 
We can use $f$ as a \textit{polynomial possibility density function} when we have a corresponding possibilistic uncertainty random variable.
\end{definition}

Consider a polynomial degree $n$ fuzzy number $a/b/c/d$, core $[b,c]$, support $[a,d]$, the corresponding fuzzy membership function (or possibility density function) which originally comes from a polynomial function $x^n$ where $n = 1,2,3, \ldots$ has the general form as shown
in table  \ref{tbl:polyDistribution}. 
\begin{table}[htp]
  \caption{General form of a polynomial possibility density function}
    \label{tbl:polyDistribution}  
    \begin{center} \renewcommand{\arraystretch}{1.2}
    \begin{tabular}{c|c}
        n    &   Left Density Function                           \\ \hline
        odd  & $f_L(x) = \frac{1}{(b-a)^n}(x-b)^n +1$ \\
        even & $f_L(x) = -\frac{1}{(b-a)^n}(x-b)^n +1$ 
        \\ \hline
        \\
        n    &  Right Density Function   
        \\ \hline
        odd    & $f_R(x) = -\frac{1}{(d-c)^n}(x-c)^n +1$ \\
        even & $f_L(x) = -\frac{1}{(d-c)^n}(x-c)^n +1$ \\ \hline        
    \end{tabular}
    \end{center}
\end{table}

In general, the upper cumulative (Possibility measure) and the lower cumulative (Necessity measure) distribution functions generated by polynomial possibility density function are $F^u$ and $F^l$, respectively, where
\[F^u(x) = Pos(x) = \left\{
\begin{array}{ll}
0 &; x < a\\
f_L(x) &; a \leq x \leq b \\
1 &; b < x
\end{array}\right.
\]
and
\[F^l(x) = Nec(x) = \left\{
\begin{array}{ll}
0 &; x < c\\
1-f_R(x) &; c \leq x \leq d \\
1 &; d < x.
\end{array}
\right.\]
For the odd number of  $n$, we calculate the left and right density functions respectively as follows:
\[d^{-}(x) = \frac{d}{dx}\left[F^u(x)\right] = 
\left\{
\begin{array}{cl}
\frac{n}{(b-a)^n}(x-b)^{n-1}&; a \leq x \leq b\\
0 &; \mathrm{ otherwise},
\end{array}\right.
\]
\[d^{+}(x) =\frac{d}{dx}\left[F^l(x)\right] =  
\left\{
\begin{array}{cl}
\frac{n}{(d-c)^n}(x-c)^{n-1}\;\;\;; c \leq x \leq d\\
0 \;\;\;\;\;\;\;\; \;\;\;\;\;\;\;\;\;\;\;\;\;\;\;\;\;\;\;\;; \mathrm{ otherwise}.
\end{array}\right.
\]
The upper cumulative distribution produces the lower integral and the lower cumulative    distribution produces the upper integral. Therefore
\[\int_{-\infty}^{\infty}x d^{-}(x)dx 
=a+ \frac{b-a}{n+1}\]
and 
\[\int_{-\infty}^{\infty}x d^{+}(x)dx 
= d - \frac{d-c}{n+1}.\]
Similar work could be done for the even degree. Thus, the interval expected value of an IVPM constructed by a polynomial possibility density function degree $n$ where $n \in \mathrm{N}$ is
\[\int_\mathrm{R} x di_{\breve{m}} = 
\left[a+ \frac{b-a}{n+1},d - \frac{d-c}{n+1}\right].\]
\begin{remark} \label{rem:regPoss}
The interval expected value of an IVPM constructed from a constant $c$ or an interval $[a,b]$, (uniform regular possibility density function), is the constant  or the interval themselves.
\end{remark}
\begin{remark} \label{rem:bound}
For any continuous possibility density function with fuzzy number $a/b/c/d$ whose core is $[b,c]$ and support is $[a,d]$, the interval expected value $\int_{\mathrm{R}}x di_{\breve{m}} =[\alpha,\beta] \subseteq [a,d]$, where $\alpha \in [a,b] $ and $\beta \in[c,d]$.

\textit{Proof:} It is a property of an expected value. $\diamond$
\end{remark}

\section{IVPM with Linear Programming} \label{lp}
Consider a linear programming (LP) with some uncertainty coefficients
\begin{eqnarray}
&\mathrm{max}&f(\vec{x},\breve{a}):= \sum_{i=1}^n \breve{a}_{i}x_i \label{eq:LP} \\
&\mathrm{s.t.}&  g(\vec{x},\breve{b},\breve{c}):= \sum_{i=1}^n \breve{b}_{ji}x_i + \breve{c}_{ji}= 0;  j = 1,\ldots, m \nonumber\\
&& \vec{0} \leq \vec{x} \leq \vec{t}; \;\;\; \mathrm{where} \;\; \vec{t} \geq \vec{0},
\nonumber
\end{eqnarray}
where some components of $\breve{a}$, $\breve{b}$ and $\breve{c}$ could represent probabilistic, possibilistic, or interval uncertain  random variables. The bound on vector $\vec{x}$ could be $\infty$. In this paper we consider the situation when these random variables are mutually independent. 
\begin{definition} (Lodwick and Jamison \cite{LodJam2007})
The IVPM constructed from two uncertain independent random variables, $\breve{X}$ and $\breve{Y}$, is defined as 
\[i_{\breve{X} \times \breve{Y}}\left(A \times B\right)\equiv i_{\breve{X}}(A) i_{\breve{Y}}(B).\]
\end{definition}

 Unlike in the ordinary LP, for this paper a problem might have no feasible region. Instead, the constraints $g(\vec{x},\breve{b},\breve{c}) = 0$ mean that $g$ can come as close to zero as possible. Therefore, it is reasonable to use penalty strategies for this type of problem. For example in \cite{LodJam2007}
\begin{eqnarray*}
&\mathrm{max}&f(\vec{x},\breve{a}):= 8x_1 +7x_2\\
&\mathrm{s.t.}& g_1(\vec{x},\breve{b},\breve{c}):= 3x_1+[1,3]x_2 + 4 = 0\\
&&g_2(\vec{x},\breve{b},\breve{c}):= \hat{2}x_1+5x_2 +1  = 0\\
&& \vec{x} \in [0,2]\;,
\end{eqnarray*}
where $\hat{2} = 1/2/2/3$. It is easy to see that this problem has no feasible set.  Therefore, the solution $\vec{x}^{*}$ for a modified problem does not need to satisfy the constraints of the original one.

Lodwick and Jamison \cite{LodJam2007} present an idea to deal with problem (\ref{eq:LP}) which involves interval types of uncertainties by using the IVPM and the operations in the order below. The explanation of these steps follows:
\begin{enumerate}
\item apply a penalty cost $\vec{p} > 0 $ (determined by a decision maker) to $\left|g(\vec{x},\breve{b},\breve{c})\right|$ ,
\item calculate the interval expected value \\$\int_{\mathrm{R}} h\left(\vec{x},\breve{a},\breve{b},\breve{c}\right)di_{\breve{a}\times\breve{b}\times\breve{c}}\;$ ( combining all liked terms $x_i$'s together or not is depending on the user), where $\;h\left(\vec{x},\breve{a},\breve{b},\breve{c}\right) = f(\vec{x},\breve{a})- \vec{p} \cdot \left|g(\vec{x},\breve{b},\breve{c})\right|,$
\item apply the ordered function $v_{[\alpha,\beta]}: [\alpha,\beta]\rightarrow \mathrm{R}^k$ to the interval coefficient $[\alpha,\beta]$ of $\int_{\mathrm{R}} h\left(\vec{x},\breve{a},\breve{b},\breve{c}\right)di_{\breve{a}\times\breve{b}\times\breve{c}}\;$ from step 2,
\item determine $u_{[\alpha,\beta]}\left(v_{[\alpha,\beta]}\left(\left[\alpha,\beta\right]\right)\right)$ where $\;u_{[\alpha,\beta]}:\mathrm{R}^k \rightarrow \mathrm{R}$,
\end{enumerate}
By assuming penalty cost vector $\vec{p}> 0$
to the constraints $\left|g(\vec{x}, \breve{b},\breve{c})\right|$, the problem becomes 
\begin{eqnarray}
\max \;\;\; h\left(\vec{x},\breve{a},\breve{b},\breve{c}\right) := f(\vec{x},\breve{a}) - \vec{p}\cdot \left|g(\vec{x},\breve{b},\breve{c})\right| \label{eq:h_abs}
\end{eqnarray}
For now, we ignore the fact that $\left|g(\vec{x},\breve{b},\breve{c})\right|$ is a non smooth function and consider it as a linear function, so that we can explain clearly each step of the operations above.
Next, by using the assumption that all uncertain random variables are independent, we form an IVPM $i_{\breve{a}\times\breve{b} \times \breve{c}}$ then calculate the interval expected value with respect to this IVPM. We get
\begin{eqnarray}
\int_{\mathrm{R}}  h\left(\vec{x},\breve{a},\breve{b},\breve{c}\right) di_{\breve{a}\times \breve{b} \times \breve{c}}\;.
\label{eq:interval-expected value}
\end{eqnarray}
Since $i_{\breve{a}\times\breve{b}\times \breve{c}}\equiv i_{\breve{a}}\; i_{\breve{b}}\;i_{\breve{c}}$ ($\breve{a}$, $\breve{b}$ and $\breve{c}$ are independent), the equation (\ref{eq:interval-expected value}) is a triple integral with respect to $i_{\breve{a}}$, $i_{\breve{b}}$ and $i_{\breve{c}}$. We consider $x_i$ as a constant while finding the interval expected value for the IVPM. We get a linear function with interval coefficients as a result for this step.  

Noting that if we  find the interval expected value for each function $f$ and $g$ before applying the violation cost, the result up to this point is the same as working with violation cost then the interval expected value of function $h$.
So far we achieve a linear unconstrained objective function with interval coefficients.

The next question is how to deal with these coefficients to keep as much information about the interval as possible. This information is dependent on the decision maker. He might use the midpoint as his priority, or he might want to keep track on the length of the interval. So the decision maker needs to put his priorities in the order. For example, he might use the midpoint as his first priority since it is the best estimate for the true value. Then his second priority could be the length of this interval because together with the midpoint, he will be able to get back to his interval easily. His third priority could be the right end point of the interval because he did not want to excess that limit, and so on. The decision maker also can have different orders and (or) methods for each of intervals.

The decision maker will write down the function $v_{[\alpha,\beta]}:[\alpha,\beta]\rightarrow \mathrm{R}^k$  to represent his $k$ priorities. For example when $k = 3$,
\[v_{[\alpha,\beta]}\left(\left[\alpha,\beta\right]\right) = \left(
\begin{array}{c}
\frac{\alpha+\beta}{2}\\
\beta-\alpha\\
\beta
\end{array}
\right).\]
For the different interval coefficient $[\mu,\nu]$, the corresponding function $v_{[\mu,\nu]}$ might have different orders from $v_{[\alpha,\beta]}$.

Now, the decision maker might weigh these priorities equally or might have some fancy strategy to deal with them. Again, these strategies depend on the interval $[\alpha,\beta]$ and the opinion of the decision maker. For example, given equal importance, he can define the function $u_{[\alpha,\beta]}: \mathrm{R}^k \rightarrow \mathrm{R}$ as 
\[u_{[\alpha,\beta]}\left(\begin{array}{c}
\frac{\alpha+\beta}{2}\\
\beta-\alpha\\
\beta,
\end{array}
\right) = \frac{1}{3}\left(\frac{\alpha+\beta}{2}\right)+\frac{1}{3}(\beta-\alpha)+\frac{1}{3}\beta.\] 
 The bottom line is that now the interval coefficients become a real number by using the functions $u_{[\alpha,\beta]} $ and $v_{[\alpha,\beta]}$ on  the interval coefficient $[\alpha,\beta]$, i.e. \[u_{[\alpha,\beta]}\left(v_{[\alpha,\beta]}\left(\left[\alpha,\beta\right]\right)\right) \in \mathrm{R}.\]

Therefore, we transform our original problem to a linear unconstrained objective function problem with real coefficients. Together with the bound on $\vec{x}$ we get a solution for our transformation problem, (drop the subscript $[\alpha,\beta]$ from functions $u$ and $v$),
\[\max_{x} \;\;u\left(v\left(\int_{\mathrm{R}}h(\vec{x},\breve{a},\breve{b}, \breve{c})\;di_{\breve{a}\times \breve{b}\times \breve{c}}\right)\right).\]
Unfortunately, we have to deal with the non smooth function,  $h(\vec{x},\breve{a},\breve{b}, \breve{c})$. Also we will not be that lucky to get a linear unconstrained objective function after these operations. 

By looking at step $2$ suggested above carefully, we can calculate the interval expected value of each constraint $g_i$ and objective function $f$ before applying a penalty cost $\vec{p}$ to them (without changing the result at the end of step 4). So, there will be some changes in the operations stated above. The operations we use in this paper are in the following steps.
\begin{algorithm} \label{alg:LP} IVPM with interval uncertainty coefficients LP. 
\begin{enumerate}
\item Calculate the interval expected value of each $g_i(\vec{x},\breve{a}, \breve{b},\breve{c})$, and $f(\vec{x},\breve{a})$, i.e., find $\int_{\mathrm{R}}g_i(\vec{x}, \breve{b},\breve{c}) \;di_{\breve{a}\times \breve{b}\times \breve{c}}$ for $i = 1, \ldots, m$ and $\int_{\mathrm{R}}f(\vec{x},\breve{a})\;di_{\breve{a}\times \breve{b}\times \breve{c}}$\;. For convenience, we store the result of this step as
\begin{eqnarray*}
f(\vec{x}) &\leftarrow& \int_{\mathrm{R}}f(\vec{x},\breve{a})\;di_{\breve{a}\times \breve{b}\times \breve{c}}\\
g_i(\vec{x}) &\leftarrow& \int_{\mathrm{R}}g_i(\vec{x}, \breve{b},\breve{c})\;di_{\breve{a}\times \breve{b}\times \breve{c}}\;; \;\;\; i = 1,\ldots,m\;.
\end{eqnarray*}
\item Apply the ordered function $v_{[\alpha,\beta]}: [\alpha,\beta]\rightarrow
\mathrm{R}^k $ to the interval coefficient $[\alpha,\beta]$ of $f(\vec{x})$ and $g_i(\vec{x})$ received from step 1. Again, store the result as
\begin{eqnarray*}
f(\vec{x}) &\leftarrow& v\left(\int_{\mathrm{R}}f(\vec{x},\breve{a})\;di_{\breve{a}\times \breve{b}\times \breve{c}}\right)\\
g_i(\vec{x}) &\leftarrow& v\left(\int_{\mathrm{R}}g_i(\vec{x}, \breve{b},\breve{c})\;di_{\breve{a}\times \breve{b}\times \breve{c}} \right);  i = 1,\ldots,m.
\end{eqnarray*}
\item Determine $u_{[\alpha,\beta]}\left(v_{[\alpha,\beta]}\left([\alpha,\beta]\right)\right) \in \mathrm{R}$ to the coefficient $v_{[\alpha,\beta]}\left([\alpha,\beta]\right)$ of $f(\vec{x})$ and $g_i(\vec{x})$ received from step 2, where $u_{[\alpha,\beta]}:\mathrm{R}^k \rightarrow \mathrm{R}$.
Store the result as
\begin{eqnarray*}
f(\vec{x}) &\leftarrow& u\left(v\left(\int_{\mathrm{R}}f(\vec{x},\breve{a})\;di_{\breve{a}\times \breve{b}\times \breve{c}}\right)\right)\\
g_i(\vec{x}) &\leftarrow& u\left(v\left(\int_{\mathrm{R}}g_i(\vec{x}, \breve{b},\breve{c})\;di_{\breve{a}\times \breve{b}\times \breve{c}}\right)\right). 
\end{eqnarray*}
Now, the coefficients of $f$ and $g_i$ are all constants and we approach the following LP problem
\begin{eqnarray}
\max \;\;&f(\vec{x})  \nonumber\\
\mathrm{s.t.} \;\; &g_i(\vec{x}) = 0;\;\;\; i = 1, \ldots , m  \label{eq:nopenalty}\\
	& \vec{0} \leq \vec{x} \leq \vec{t}; \;\;\; \mathrm{where} \;\; \vec{t} \geq \vec{0},
\nonumber				    
\end{eqnarray}

\item Apply a penalty cost $\vec{p} > 0$ to the function vector $g$ received from step 3. Define the penalty function $h(\vec{x})= f(\vec{x}) - p^T\left|g(\vec{x})\right|$. The problem becomes an unconstrained objective function to maximize the function $h(\vec{x})$.
\item Use modeling techniques to get rid of non-smooth absolute function, $\left|g(\vec{x})\right|$.
\end{enumerate}
\end{algorithm}
In general, the penalty vector $\vec{p}$ is not a fixed vector. It depends upon the excess or shortage of function $g$ from $0$. Moreover, one of the modeling difficulty is that we could not know in advance which constraint will lack or excess the balance zero. Let us denote \$$e_i$ and \$$s_i$ as the cost penalty for each excess and shortage unit of $g_i(\vec{x})$ from zero, respectively. Noting that $\left|g(\vec{x},\breve{b},\breve{c})\right|$ , the absolute of the function vector $g$, is a non smooth function. So we define $\psi_i :=  \max\;\left\{0, g_i(\vec{x})\right\}$ and 
$\zeta_i := \max\;\left\{0, -g_i(\vec{x})\right\}$, for $i = 1,\ldots, m$. Then the unconstrained problem got form step 4 of Algorithm \ref{alg:LP}
\[\max \; h(\vec{x}):= f(\vec{x}) - p^T \left|g(\vec{x})\right|\]
can be remodeled as
\begin{eqnarray}
\max \;& h(\vec{x},\vec{\psi},\vec{\zeta}) := f(\vec{x}) - \sum_{i=1}^m e_i \psi_i - \sum_{i=1}^m s_i \zeta_i \nonumber\\
 \label{eq:model_lp}\\
\mathrm{s.t.} &\left.\begin{array} {l}
\left.\begin{array} {l}
              \psi_i \geq g_i(\vec{x}) \\
					     \zeta_i \geq -g_i(\vec{x})\\
					     \psi_i \geq 0 \\
					     \zeta_i \geq 0 
					    \end{array} \right\} \;\; i = 1,\ldots,m  ;\\
					     \\
					    0 \leq \vec{x} \leq \vec{t};\;\;\; \vec{t} \geq 0 .\end{array}\right.\nonumber
					  					    \end{eqnarray}
From the explanation above we provide the conclusion in the following lemma.
\begin{lemma} \label{lem:LP}
Consider an LP (\ref{eq:LP}), after working through the algorithm \ref{alg:LP},
the remodeled problem (\ref{eq:model_lp}) becomes an LP problem.

\textit{Proof:} It is clear from the algorithm \ref{alg:LP}.
\end{lemma}
\begin{example} \label{ex:lp} Consider the problem
\begin{eqnarray*}
\max \;& f(\vec{x},\breve{a}) &:= \widehat{2}x_1 - \bar{3}x_2+\left[3,5\right]x_3\\
\mathrm{s.t.} \;& g_1(\vec{x},\breve{b},\breve{c})&:=\widehat{4}x_1 + \left[1,5\right]x_2-2x_3-[0,2] =0\\
& g_2(\vec{x},\breve{b},\breve{c})&:= 6x_1 -\bar{2}x_2 +9x_3 -9 = 0\\
&g_3(\vec{x},\breve{b},\breve{c})&:=-2x_1-[1,4]x_2-\widehat{8}x_3 +\bar{5}=0\\
&& 0 \leq x_1 \leq 3\\
&& 0 \leq x_2 \leq 2 \\
&& 0 \leq x_3 \leq 2,
\end{eqnarray*}
where $\breve{a} = \left[\widehat{2},-\bar{3},\left[3,5\right]\right]^T$, 
$\breve{b} = \left[\begin{array} {ccc}
  \widehat{4} & \left[1,5\right] & -2 \\
  6 & -\bar{2}& 9\\
  -2 & -[1,4] & -\widehat{8}
\end{array} \right]$ and  $\breve{c} =\left[-\left[0,2\right],-9,\bar{5}\right]^T$.\\
\textit{Note}: The functions $f$ and $g_3$ involve 3 types of uncertainties.

These coefficients have possibility (or probability) polynomial density functions as shown in table \ref{tbl:exlp}. Apply Algorithm \ref{alg:LP} step by step, we have:
\begin{table}[htp]
    \caption{ Polynomial possibility density function for the coefficients}
    \label{tbl:exlp}  
    \begin{center} \renewcommand{\arraystretch}{1.35}
    \begin{tabular}{c|c|c}
        coefficient   &  a/b/c/d   &  degree of polynomial distribution   \\ \hline
        $\widehat{2}$  & $0/1/2/3$ & $2$ \\\hline
        $\bar{3}$ & $2/3/3/4$& $1$ \\ \hline
        $\widehat{4}$ & $ 2/4/4/6$ & $1$ \\ \hline
        $\bar{2}$ & $1/2/2/3$ & $1$ \\ \hline
        $\widehat{8}$ &$7/8/8/9$& $3$ \\\hline
        $\bar{5}$ &$4/5/5/6$&$1$ \\ \hline         
    \end{tabular}
    \end{center}
\end{table}

\vspace{-.28cm}
\begin{enumerate}
\item  Calculate the interval expected values
\begin{eqnarray*}
f(\vec{x}) &\leftarrow& 
\left[\frac{1}{3},\frac{8}{3}\right] x_1 - \left[3,3\right] x_2 + \left[3,5\right]x_3\\
&=& \left[\frac{1}{3},\frac{8}{3}\right] x_1 + \left[-3,-3\right] x_2 + \left[3,5\right]x_3 \; ;\\
g_1(\vec{x}) &\leftarrow& 
\left[3,5\right] x_1 +\left[1,5\right]x_2 -2x_3  - [0,2]\\
&=& \left[3,5\right] x_1 +\left[1,5\right]x_2 +[-2,-2]x_3 + \\ &&[-2, 0] \; ;\\
g_2(\vec{x}) &\leftarrow& 
6 x_1 -\left[2,2\right]x_2 +9x_3  -9\\
&=& [6,6] x_1 +\left[-2,-2\right]x_2 +[9,9]x_3 \\ && +[-9,-9] \; ;\\
g_3(\vec{x}) &\leftarrow& 
 -2x_1 -\left[1,4\right]x_2 -\left[\frac{29}{4},\frac{35}{4}\right]x_3  +5\\
&=&  [-2,-2]x_1 +\left[-4,-1\right]x_2 + \\ &&\left[-\frac{35}{4},-\frac{29}{4}\right]x_3  +[5,5].
\end{eqnarray*}
\item  Assume that $v\left[\alpha,\beta\right] = \left(\frac{\alpha+\beta}{2}, \beta-\alpha\right)$, therefore 
\begin{eqnarray*}
f(\vec{x}) &\leftarrow& 
(\frac{3}{2},\frac{7}{3})x_1 +(-3,0)x_2+(4,1)x_3\;;\\
g_1(\vec{x}) &\leftarrow& 
(4,1)x_1 +(3,4)x_2+(-2,0)x_3 \\&&+(-1,2)\;;\\
g_2(\vec{x}) &\leftarrow& 
(6,0)x_1 +(-2,0)x_2+(9,0)x_3\\&&+(-9,0) \;;
\end{eqnarray*}
\begin{eqnarray*}
g_3(\vec{x}) &\leftarrow&
 (-2,0)x_1+(-\frac{5}{2},3)x_2 +(-16,\frac{3}{2})x_3\\ &&+(5,0).
\end{eqnarray*}
\item  The user applies the function $u:\mathrm{R}^2 \rightarrow \mathrm{R}$ as $u \left(\left(a_1, a_2\right)\right) = \frac{a_1+a_2}{2}$, when the original constants are uncertain and $u \left(\left(a_1, a_2\right)\right) = a_1$, when the original coefficients are constants.
\begin{eqnarray*}
f(\vec{x}) &\leftarrow& 
\frac{23}{6}x_1 -\frac{3}{2}x_2+\frac{5}{2}x_3\;;\\
g_1(\vec{x}) &\leftarrow& 
 \frac{5}{2}x_1 +\frac{7}{2}x_2-2x_3 +\frac{1}{2}\;;\\
g_2(\vec{x}) &\leftarrow& 
 6x_1 -x_2+9x_3-9 \;;\\
g_3(\vec{x}) &\leftarrow& 
 -2x_1+\frac{1}{4}x_2 -\frac{29}{4}x_3+\frac{5}{2}.
\end{eqnarray*}
\item  Choose $s_1 = s_2= s_3 =$ \$2, and $e_1 =e_2=e_3 =$ \$1. The problem becomes
\[\max\; h(\vec{x}):=  \frac{23}{6}x_1 -\frac{3}{2}x_2+\frac{5}{2}x_3 -\vec{p}^T \left|g(x)\right|.\]
\item  Remodel the problem in step 4:
\begin{eqnarray*}
\max \; &h(\vec{x},\vec{\psi},\vec{\zeta}) &:=\; \frac{23}{6}x_1 -\frac{3}{2}x_2+\frac{5}{2}x_3 -2\psi_1\\
&&\;\;\;\;\;-2\psi_2-2\psi_3-\zeta_1-\zeta_2-\zeta_3\\
\mathrm{s.t.} &\psi_1 &\geq \;  \frac{5}{2}x_1 +\frac{7}{2}x_2-2x_3 +\frac{1}{2}\\
& \psi_2& \geq \;   6x_1 -x_2+9x_3-9 \\
& \psi_3& \geq \;  -2x_1+\frac{1}{4}x_2 -\frac{29}{4}x_3+\frac{5}{2}\\
&\zeta_1& \geq \;  -\frac{5}{2}x_1 -\frac{7}{2}x_2+2x_3 -\frac{1}{2}\\
&\zeta_2& \geq \;  -6x_1 +x_2-9x_3+9\\
&\zeta_3& \geq \;  2x_1-\frac{1}{4}x_2 +\frac{29}{4}x_3-\frac{5}{2}\\
&& 0 \leq x_1 \leq 3\\
&& 0 \leq x_2 \leq 2 \\
&& 0 \leq x_3 \leq 2\\
&& \psi_i \geq 0, \;\; \zeta_i \geq 0;\;\;\;\;\; i = 1,2,3.
\end{eqnarray*} 
\end{enumerate}
Solving this problem using GAMS, we have an optimal solution at\\ $\left(x_1,x_2,x_3,\psi_1,\psi_2,\psi_3,\zeta_1,\zeta_2,\zeta_3\right) = \left(0.3913,0,0.7391,0,0,0,0,0,3.6413\right)$
and the optimal objective value is -0.2935.
\end{example}
\section{IVPM with Quadratic Programming} \label{nlp}
In this section we start with a quadratic programming (QP) with some independent uncertainty coefficients as follow, 
\begin{eqnarray}
&\mathrm{max}&f(\vec{x},\breve{M},\breve{b}):= \frac{1}{2}\vec{x}^T \breve{M}\vec{x} + \breve{b}^T \vec{x}\label{eq:QP}\\
&\mathrm{s.t.}& g(\vec{x},\breve{c},\breve{d}):= \sum_{i=1}^n \breve{c}_{ji}x_i + \breve{d}_{ji}= 0; j = 1, \ldots, m \nonumber\\
&& \vec{0} \leq \vec{x} \leq \vec{t}; \;\;\; \mathrm{where} \;\; \vec{t} \geq \vec{0}.
\nonumber
\end{eqnarray}
Note that 
\[\breve{M}= 
\left[\begin{array}{ccc}
\breve{a}_{11}& \cdots & \breve{a}_{1n}\\
\vdots & \ddots &\vdots\\
\breve{a}_{n1}& \cdots & \breve{a}_{nn}
\end{array}
\right]\]
is not necessarily a symmetric matrix. Also, $f(\vec{x},\breve{M},\breve{b})$ is a quadratic function with uncertain coefficients and $ g(\vec{x},\breve{c},\breve{d})$ is a linear function vector.
 Elements of $\breve{M}$, ${\breve{b}}$, ${\breve{c}}$ and ${\breve{d}}$ could be constant, possibilistic, fuzzy or interval uncertain random variables. 
 
By applying the Algorithm \ref{alg:LP} to this QP, the remodeled problem  for (\ref{eq:QP}) is similar to the problem (\ref{eq:model_lp}). The only difference is that the function $f$ 
is now a quadratic function. So we can rewrite  the problem (\ref{eq:model_lp}) as
\begin{eqnarray}
\max \;& h(\vec{x},\vec{\psi},\vec{\zeta}) := \frac{1}{2}\vec{x}^T M\vec{x} + \vec{b}^T\vec{x} - \sum_{i=1}^m e_i \psi_i\nonumber\\ & - \sum_{i=1}^m s_i \zeta_i \nonumber\\
 \label{eq:model_qp}\\
\mathrm{s.t.} &\left.\begin{array} {l}
\left.\begin{array} {l}
              \psi_i \geq g_i(\vec{x}) \\
					     \zeta_i \geq -g_i(\vec{x})\\
					     \psi_i \geq 0 \\
					     \zeta_i \geq 0 
					    \end{array} \right\} \;\; i = 1,\ldots,m  ;\\
					     \\
					    0 \leq \vec{x} \leq \vec{t};\;\;\; \vec{t} \geq 0 .\end{array}\right.\nonumber
					  					    \end{eqnarray}
\vspace{.18cm}

where $M$ and $\vec{b}$ are generated from step 3 and step 4 of Algorithm \ref{alg:LP}. Function $g$ is now a linear constrained function vector whose coefficients are constants. Without loss of generality, $M$ can be a symmetric matrix. If M is positive or negative definite, it will be fairly easy to solve this remodeled problem (\ref{eq:model_qp}). However, we can not get a nice form of matrix $M$ in general. It depends on the types of uncertain coefficients and how the decision maker defines function $u $ and $v$. (In this case, we can use technique of elimination of variables using linear constraints, \cite{NocWri2000} , to get at least a stationary point. Another method we could consider is one using complementarity problem and modify the QP problem to a linear problem with complementarity constraints \cite{Baz_etal1993}.)

From the work  in section 3 and 4, we can conclude that if we start with an optimization with uncertain coefficients and linear constraints, the Algorithm \ref{alg:LP} leads us to the same type of optimization problem with constant coefficients. 

We illustrate an example of a QP problem by changing the objective function of example \ref{ex:lp}.
\begin{example} \label{ex:nlp} We use the example \ref{ex:lp}. The only difference is the objective function, 
$ \mathrm{i.e.,}\;\;f(\vec{x},\breve{a}) := \widehat{2}x_1 - \bar{3}x_2+\left[3,5\right]x_3\;\; \mathrm{becomes} \;\;
 f(\vec{x},\breve{a}) := \widehat{2}x_1^2 - \bar{3}x_2^2+\left[3,5\right]x_3.
$

By go through the steps in the algorithm \ref{alg:LP}, the calculation remains the same as one in example \ref{ex:lp}. Therefore, in the step 5 we achieve
\vspace{.2cm}
\begin{eqnarray*}
\max \; &h(\vec{x},\vec{\psi},\vec{\zeta}) &:=\; \frac{23}{6}x_1^2 -\frac{3}{2}x_2^2+\frac{5}{2}x_3 -2\psi_1\\
&&\;\;\;\;\;-2\psi_2-2\psi_3-\zeta_1-\zeta_2-\zeta_3\\
\mathrm{s.t.} &\psi_1 &\geq \;  \frac{5}{2}x_1 +\frac{7}{2}x_2-2x_3 +\frac{1}{2}\\
& \psi_2& \geq \;   6x_1 -x_2+9x_3-9 \\
& \psi_3& \geq \;  -2x_1+\frac{1}{4}x_2 -\frac{29}{4}x_3+\frac{5}{2}\\
&\zeta_1& \geq \;  -\frac{5}{2}x_1 -\frac{7}{2}x_2+2x_3 -\frac{1}{2}\\
&\zeta_2& \geq \;  -6x_1 +x_2-9x_3+9\\
&\zeta_3& \geq \;  2x_1-\frac{1}{4}x_2 +\frac{29}{4}x_3-\frac{5}{2}\\
&& 0 \leq x_1 \leq 3\\
&& 0 \leq x_2 \leq 2 \\
&&0 \leq x_3 \leq 2\\
&& \psi_i \geq 0, \;\; \zeta_i \geq 0;\;\;\;\;\; i = 1,2,3.
\end{eqnarray*} 
Solving this problem using GAMS, we have an optimal solution at\\ $\left(x_1,x_2,x_3,\psi_1,\psi_2,\psi_3,\zeta_1,\zeta_2,\zeta_3\right) = \left(0.3913,0,0.7391,0,0,0,0,0,3.6413\right)$
and the optimal objective value is -1.2065.
\end{example}
\section{Conclusion and Further Research}

What we have done in this paper is that we use the concept of IVPMs to get the corresponding interval expected value of  uncertain coefficients in an LP (or a QP with linear constraints) problem. Then our optimization problem becomes the same type as the original problem with interval coefficients. By using functions $u$ and $v$ (given by a decision maker), our problem is an ordinary problem which can be solved by any appropriate tools such as GAMS and LINDO. Moreover, it implied that IVPM (with independent uncertain random variables) can be used to  put all three uncertainties, (probabilistic, possibilistic and interval), in one framework. So that optimization problems containing all three types of uncertainties in one
problem, especially in one constraint inequality, (constraint $g_3$ in example \ref{ex:lp}, for instance), can be solved. It needs an extra work before we can conclude that a similar statement holds (or not) if we add the other types of uncertainties to the problems.

We restricted our uncertain random variables to be independent which makes it much easier when calculating the interval expected values. The question is that can we still be able to use the concept of IVPMs when we have dependent uncertain random variables. Therefore, the suggestion for further research is focusing on the dependence of uncertain coefficients. Also, we need to give a concrete semantics for our work.

\nocite{*}
\renewcommand{\refname}{Bibliography}

\addcontentsline{toc}{section}{\refname} 
\bibliography{strings,NLP}      
\bibliographystyle{plainPlus}       
\end{document}